\font\smc=cmcsc10

\def \varpi {\bar \omega}
\def \N {I \! \! N}
\def \Z{Z \! \! \! Z}

\def \N {I \! \! N}

\def\text#1{{\textstyle\rm#1}}

\def\newline{\unskip\null\hfill\break}
\def\square{\sqcup\!\!\!\!\sqcap}

\def\cejour{
\number\day\space\ifcase\month\or janvier\or f\'evrier \or
mars \or avril\or mai \or juin\or juillet\or ao\^ut\or
septembre\or octobre\or novembre\or
d\'ecembre\fi\space\number\year} 
\def\today{
\space\ifcase\month\or January\or February \or
March \or April\or May \or June\or July\or August\or
September\or October\or November\or
December\fi\space
\number\day ,\space\number\year} \nonstopmode
\def\cejour{(\number\time)
\number\day\space\ifcase\month\or janvier\or f\'evrier \or
mars \or avril\or mai \or juin\or juillet\or ao\^ut\or
septembre\or octobre\or novembre\or
d\'ecembre\fi\space\number\year} 
\def\C{{\mathchoice {\setbox0=\hbox{$\displaystyle\rm C$}\hbox{\hbox
to0pt{\kern0.4\wd0\vrule height0.9\ht0\hss}\box0}}
{\setbox0=\hbox{$\textstyle\rm C$}\hbox{\hbox
to0pt{\kern0.4\wd0\vrule height0.9\ht0\hss}\box0}}
{\setbox0=\hbox{$\scriptstyle\rm C$}\hbox{\hbox
to0pt{\kern0.4\wd0\vrule height0.9\ht0\hss}\box0}}
{\setbox0=\hbox{$\scriptscriptstyle\rm C$}\hbox{\hbox
to0pt{\kern0.4\wd0\vrule height0.9\ht0\hss}\box0}}}}

\def\frac#1#2{{#1\over#2}}
\def\square{\sqcup\!\!\!\!\sqcap}
\def\sec{'\! '}


\noindent \today  \

\vskip 1mm

\centerline { {\smc About the characterization of some residue currents}}
\vskip 2mm\centerline{Pierre Dolbeault}\vskip 5mm

This unpublished paper is a copy (completed by a development of section 5 and by minor corrections) of  the article with the same title published in:

Complex Analysis and Digital Geometry, {\it  Proceedings from the Kiselmanfest, 2006}, Acta Universitatis Upsaliensis, C. Organisation och Historia, {\bf 86}, Uppsala University Library (2009), 147-157.\vskip 3mm

\centerline {\smc Contents}\vskip 3mm

1. Introduction

2. Preliminaries: local description of a residue current

3. The case of simple poles

4. Expression of the residue current of a closed differential form

5. Generalization of Picard's theorem. Structure of residue currents of closed meromorphic forms

6. Remarks about residual currents

\vskip 5mm

\noindent {\bf 1. Introduction}.

\vskip 2mm

\noindent {\bf 1.1. Residue current in dimension 1.}
Let $\omega= g(z) dz$ be a meromorphic 1-form on a small enough 
 open set $0\in U\subset \C$ having 0 as unique pole, with multiplicity $k$:
 
 $$g=\sum_{l=1}^k \frac {a_{-l}}{z^l}+{\rm holomorphic  \ function}$$

Note that $\omega$ is $d$-closed.\vskip 1mm
 
 Let $\psi=\psi_0 d\overline z\in {\cal D}^1(U)$ be a 1-test form. In general $g\psi$ is not integrable, but the principal value
 
 $$Vp\lbrack\omega\rbrack(\psi) =\lim_{\epsilon\rightarrow 0}\int_{\vert z\vert\geq\epsilon}\omega\wedge\psi \ \  
$$
 
 exists, and $d Vp\lbrack\omega\rbrack=d\sec Vp\lbrack\omega\rbrack={\rm Res}\lbrack\omega\rbrack$ is the residue current of 
$\omega$. For any test function $\varphi$ on $U$,
 
 $${\rm Res}\lbrack\omega\rbrack(\varphi)=\lim_{\epsilon\rightarrow 0}\int_{\vert z\vert=\epsilon}\omega
\wedge\varphi$$
 
 Then ${\rm Res}\lbrack\omega\rbrack=2\pi i \ {\rm res}_0(\omega)\delta_0+
dB=\displaystyle\sum_{j=0}^{k-1}b_j\frac{\partial^j}{\partial z^j}\delta_0$ \  where $ {\rm res}_0(\omega)= a_{-1}$
is the Cauchy residue. We remark that $\delta_0$ is the integration current on the subvariety $\{0\}$ of $U$, that
$D=\displaystyle\sum_{j=0}^{k-1}b_j\frac{\partial^j}{\partial z^j}$ and that $b_j=\lambda_j a_{-j}$ where the $\lambda_j$
are universal constants.
 
 Conversely, given the subvariety $\{0\}$ and the differential operator $D$, then the meromorphic differential form $\omega$ is equal to $gdz$, up to holomorphic form; hence the residue current ${\rm Res}\lbrack\omega\rbrack=D\delta_0$, can be constructed.\vskip 1mm 
 
 \noindent  {\bf 1.2. Characterization of holomorphic chains.} P. Lelong
(1957) proved that a complex analytic subvariety $V$ in  a complex analytic
manifold $X$ defines an integration current
$\varphi\mapsto \lbrack V\rbrack(\varphi)=\displaystyle\int_{{\rm Reg}
V}\varphi$ on $X$. More generally, a {\it holomorphic p-chain} is a
current $\displaystyle\sum_{l\in L} n_l \lbrack V_l\rbrack$ where
$n_l\in\Z$, $\lbrack V_l\rbrack$ is the integration current defined by an
irreductible $p$-dimensional complex analytic subvariety $V_l$, the
family $(V_l)_{l\in L}$ being locally finite.
 
 During more than twenty years, J. King $\lbrack$K 71$\rbrack$, Harvey-Shiffman $\lbrack$HS 74$\rbrack$, Shiffman $\lbrack$S 83$\rbrack$, H. Alexander $\lbrack$A 97$\rbrack$ succeeded in proving the following structure theorem: {\it Holomorphic $p$-chains on a complex manifold $X$ are exactly the rectifiable $d$-closed currents of bidimension $(p,p)$ on $X$.}

 In the case of section 1.1, Res $\lbrack\omega\rbrack$ is the holomorphic chain with complex coefficients $2\pi i \ {\rm res}_0(\omega)\delta_0$ if and only if 0 is a simple pole of $\omega$.\vskip 2mm
 
 \noindent {\bf 1.3.} Our aim is to characterize residue currents using
rectifiable currents with coefficients that are principal values of meromorphic
differential forms and holomorphic differential operators acting on them.
 
 We present a few results in this direction.
 
 The structure theorem of section 1.2 concerns complex analytic varieties
and closed currents. So, after generalities on residue currents of
semi-meromorphic differential forms, we will concentrate on residue
currents of closed meromorphic forms.
 
 \vskip 2mm
 
\noindent  {\bf 2. Preliminaries: local description of a residue
current} ($\lbrack$D 93$\rbrack$, section 6)\vskip 1mm

\noindent {\bf 2.1.} We will consider a finite number of
holomorphic functions defined on a small enough open neighborhood $U$ of 
the origin 0 of
$\C^n$, with coordinates $(z_1,\ldots, z_n)$. For convenient coordinates, any
semi-meromorphic differential form, for $U$ small enough, can be written $\displaystyle\frac{\alpha}{f}$, where $\alpha\in{\cal E}^.(U)$), $f\in {\cal  O}(U)$ and
$$f= u_j\prod_k \ {_j\rho_k^{r_k}},$$
where the ${_j\rho_k}$ are irreducible distinct Weierstrass polynomials in $z_j$ and the $r_k\in \N$ are 
independent of $j$, moreover $u_j$ is a unit at 0, i.e., for $U$ small enough, $u_j$ does not vanish on $U$. Let $B_j$ be the discriminant of the polynomial
${_j\rho}=\prod_k \ {_j\rho_k}$ and let $Y_k=Z({_j\rho_k})$; it is clear that $Y_k$ is independent of
$j$. Let $Y=\cup_kY_k$ and $Z=$ Sing $Y$.

After shrinkage of $(0\in) \ U$, the following expressions of $\displaystyle\frac {1}{f}$ are valid on
$U$: for every $j\in \lbrack 1,\ldots,n\rbrack$,
$$\displaystyle\frac {1}{f}=u_j^{-1}\sum_k\sum_{\mu=1}^{r_k}{ ^j{c_\mu^k}}\frac{1}{{_j\rho_k^\mu}}$$
where ${ ^j{c_\mu^k}}$ is a meromorphic function whose polar set, in $Y_k$, 
is contained in $Z(B_j)$. Notice that
$B_j$ is a holomorphic function of $(z_1,\ldots, \widehat{z_j},\ldots,z_n)$. In the following, for simplicity, we omit the unit $u_j^{-1}.$\vskip 1mm

\noindent {\bf 2.2.} Let $\displaystyle\omega=\frac {1}{f}$, \hskip  3mm
$Vp
\lbrack\omega\rbrack(\psi)=\displaystyle\lim_{\epsilon\rightarrow
0}\int_{\lbrack f\rbrack\geq\epsilon}\omega\wedge\psi$; \ $\psi\in {\cal
D}^{n,n}(U)$. The residue of $\omega$ is 
$${\rm Res}\lbrack\omega\rbrack= (dVp-Vp d) \lbrack\omega\rbrack =(d\sec Vp-Vp d\sec) \lbrack\omega\rbrack$$ 

For every $\varphi\in{\cal D}^{n,n-1}(U)$, let
$\varphi=\displaystyle\sum_{j=1}^n\varphi_j$ with 
$$\varphi_j=\psi_j dz_1\wedge \ldots\wedge d\overline z_1\wedge \ldots\wedge \widehat
{d\overline z_j}\wedge\ldots$$
Then, from Herrera-Lieberman $\lbrack$HL 71$\rbrack$, and the next 
lemma about $B_j$ ,\vskip 1mm

\noindent we have: 

$${\rm Res}\lbrack\omega\rbrack(\varphi)= \sum_{j=1}^n \sum_k\sum_{\mu=1}^
{r_k}\lim_{\delta\rightarrow
0}\lim_{\epsilon\rightarrow 0}\int_{\mid B_j\mid\geq\delta 
\mid{_j\rho_k}\mid=\epsilon}{ 
^j{c_\mu^k}}\frac{1}{{_j\rho_k^\mu}}\varphi_j.$$

\noindent The lemma we have used here is the following:

\noindent {\bf Lemma 2.1. } ($\lbrack$D 93$\rbrack$, Lemma 6.2.2). 
 
$${\rm Res} \lbrack\omega\rbrack(\varphi_j)= \lim_{\delta\rightarrow
0}\lim_{\epsilon\rightarrow 0}\int_{\mid B_j\mid\geq\delta \mid 
f\mid=\epsilon}\omega\varphi_j.$$

\noindent Outside  $Z(B_j)$, for $\mid {_j\rho_k}\mid$ small enough (since
$\displaystyle\frac{\partial_j\rho_k}{\partial z_j}\not = 0$), we take $(z_1,\ldots,z_{j-1},_j\rho_k,
z_{j+1},\ldots,z_n)$ as local coordinates.\vskip 1mm

\noindent {\bf 2.3. Notations.}
For the sake of simplicity, until the end of this section, we assume $j=1$ and write $\rho_k$,
$c_\mu^k$ instead of ${_1\rho_k}$, ${^1{c_\mu^k}}$. Outside  $Z(B_1)$,  we 
take
$(\rho_k,z_2,\ldots,z_n)$ as local coordinates; then, for every $C^\infty$ function $h$ and every
$s\in \N$, we have
$$\frac{\partial^s h}{\partial\rho_k^s}=\frac{1}{(\frac{\partial\rho_k}
{\partial z_1})^{2s-1}} D_sh,\ \  {\rm for} \ s\geq 1,$$
where $D_s=\displaystyle\sum_{\alpha=1}^s\beta_\alpha^s\frac{\partial^\alpha}{\partial z_1^\alpha}$,
 $\beta_\alpha^s$ is a holomorphic function determined by $\rho_k$ and 
 $\displaystyle D_0= 
\Big(\frac{\partial\rho_k}{\partial z_1}\Big)^{-1}$.

Let

$$g_l^\mu=\displaystyle\pmatrix{
\mu-1\cr
l\cr
} \frac{1}{\Big ({\frac{\partial\rho_k}{\partial
z_1}\Big)^{2\mu-4}}} D_l\Big(\frac{c^k_\mu}{\frac{\partial\rho_k}{\partial z_1}}\Big)
,(0\leq l\leq
\mu-2);$$
$$g_{\mu-1}^\mu=\displaystyle\frac{1}{\Big ({\frac{\partial\rho_k}{\partial
z_1}\Big)^{2\mu-3}}} D_{\mu-1}\Big(\frac{c^k_\mu}{\frac{\partial\rho_k}{\partial z_1}}\Big)$$ 

Let $Vp_{Y_k,B_1}^1\lbrack g_l^\mu\rbrack$ also denote the direct image, by the inclusion $Y_k\rightarrow U$,
of the Cauchy principal value \vskip 1mm

\noindent $Vp_{Y_k,B_1}\lbrack g_l^\mu\rbrack$ of $g_l^\mu\mid_{Y_k}$;

$D_{1,k}^{\mu,l}=\displaystyle\sum_{\alpha=1}^{\mu-1-l}(-1)^\alpha\beta_\alpha^{\mu-1-l}
\frac{\partial^\alpha}{\partial z_1^\alpha}$, and $D_{1,k}^{\mu,\mu-1}=$ id.

\noindent {\bf 2.4. Final expression of the residue}.
All what has been done for $j=1$ is valid for any $j\in\{1,\ldots,n\}$: the principal value
$Vp^j(k,\mu,l)=Vp\ ^j_{Y_k,B_j}\lbrack g_l^\mu\rbrack$ defined on $Y_k$ and the holomorphic
differential operator
$D_{j,k}^{\mu,l}$. We also denote $Vp^j(k,\mu,l)$ the direct image of the principal value by the canonical injection $Y\hookrightarrow\ U$. Then, denoting L the inner product, we have: 
$${\rm Res}\lbrack\omega\rbrack(\varphi)=2\pi
i\sum_{j=1}^n\Bigl[\sum_k\sum_{\mu=1}^{r_k}\frac{1}{(\mu-1)!}\sum_{l=0}^{\mu-1}
D_{j,k}^{\mu,l} Vp^j(k,\mu,l)\Bigr]\big(\frac{\partial}{\partial z_j}{\rm L}\varphi_j\big)\leqno (*)$$\vskip 2mm

\noindent {\bf 3. The case of simple poles.}\vskip 1mm

\noindent {\bf 3.1. The case}  $\omega=\displaystyle\frac{1}{f}$.

\noindent {\bf Lemma 3.1.} {\it For a simple pole and for every
$k$, $^j{c_1^k}$ is holomorphic}.\vskip 1mm

\noindent {\it Proof}. Let $w=z_j$ and $y=(z_1,\ldots,\hat
z_j,\ldots,z_n)$. At points $z\in U$ where $B_j(z)\not =0$, for given
$y$,  let $w_{ks}, s=1,\ldots,s_k$, be the zeros of $\rho_k$. For given
$y$, $\rho_k=\displaystyle\prod_{s=1}^{s_k}(w-w_{ks})$, 
$$\displaystyle\frac {1}{f}=u_j\sum_k\sum_{s=1}^{s_k}{^j{\cal C}}_1^{k,s}
(w-w_{ks})^{-1}$$
where ${^j{\cal C}}_1^{k,s}=\displaystyle\frac{1}{\frac{\partial}{\partial 
w}f(w_{ks},y)}$; let $\prod_\sigma^s$ denote the product for all
$\sigma\not = s$,
 $$\sum_{s=1}^{s_k}{^j{\cal C}}_1^{k,s}(w-w_{ks})^{-1}=\sum_{s=1}^{s_k}{^j
{\cal C}}_1^{k,s}\frac
{\prod_\sigma^s(w-w_{k\sigma})}{\prod_\sigma(w-w_{k\sigma})}={^jc}
_1^k(w,y)\rho_k^{-1},$$
 with

${^jc}_1^k(w,y)=\displaystyle\sum_{s=1}^{s_k}\frac{\prod_\sigma^s
(w-w_{k\sigma})}{\frac{\partial}{\partial
w}f(w_{ks},y)}$ (\noindent $\lbrack$D 57$\rbrack$, IV.B.3 et C.1).
 
\noindent Here ${^jc}_1^k(w,y)$ holomorphically extends to points of $U$
where the $w_s$ are not all distinct because: if $w_s$ appears $m$ times
in $\prod_\sigma(w-w_{k\sigma})$, it appears $(m-1)$ times in the
numerator and the denominator of
$\displaystyle\frac{\prod_\sigma^s(w-w_{k\sigma})}{\frac{\partial}{\partial
w}f(w_{ks},y)}$.\hskip 1mm $\square$\vskip 1mm

\noindent All the poles of $\omega$ are simple, i.e. for every $k$, 
$r_k=1$; then $\mu=1$, $l=0$.

$${\rm Res}[\omega\rbrack(\varphi)=2\pi
i\sum_{j=1}^n\Bigl[\sum_k
D_{j,k}^{1,0} Vp^j(k,1,0)\Bigr]\big(\frac{\partial}{\partial z_j}{\rm L}\varphi_j\big)$$

But $D_{1,k}^{1,0}=$ id; \ $\displaystyle D_0= 
\Big(\frac{\partial\rho_k}{\partial z_1}\Big)^{-1}$;\  $g_{\mu-1}^\mu=\displaystyle\frac{1}{\Big
({\frac{\partial\rho_k}{\partial z_1}\Big)^{2\mu-3}}}
D_{\mu-1}\Big(\frac{c^k_\mu}{\frac{\partial\rho_k}{\partial z_1}}\Big)$;\ 
$g_0^1=\displaystyle\frac{1}{\Big ({\frac{\partial\rho_k}{\partial
z_1}\Big)^{-1}}} D_0\Big(\frac{c^k_1}{\frac{\partial\rho_k}{\partial
z_1}}\Big)$

$=\displaystyle\frac{1}{\Big ({\frac{\partial\rho_k}{\partial z_1}
\Big)^{-1}}}
\Big(\frac{\partial\rho_k}{\partial z_1}\Big)^{-1}\Big(\frac{c^k_1}
{\frac{\partial\rho_k}{\partial
z_1}}\Big)$=$\displaystyle\Big(\frac{\partial\rho_k}{\partial z_1}
\Big)^{-1} \ c^k_1$;\hskip 5mm
 
 $$Vp^j(k,1,0)=Vp\ ^j_{Y_k,B_j}[g_0^1\rbrack\displaystyle =
Vp\ ^j_{Y_k,B_j}\Bigl[
\Big(\frac{\partial\rho_k}{\partial z_j}\Big)^{-1} \ { ^j{c_1^k}}
\Bigr],$$

\noindent hence

$${\rm Res}[\omega\rbrack(\varphi)=2\pi
i\sum_{j=1}^n\Bigl[\sum_k
 Vp\ ^j_{Y_k,B_j}\lbrack
\Big(\frac{\partial\rho_k}{\partial z_j}\Big)^{-1} \ {
^j{c_1^k}}]\big(\frac{\partial}{\partial z_j}{\rm L}\varphi_j\big)
\Bigr]$$
where $^j{c_1^k}$ is holomorphic.

\noindent {\bf 3.2. The case of any degree.} Let  
$\omega=\displaystyle\frac{\alpha}{f}$.
Then Res $\lbrack \omega\rbrack =\alpha\wedge$
Res$\displaystyle(\frac{1}{f})$. Moreover, $d$ Res
$\lbrack\omega\rbrack=\pm {\rm Res}\lbrack d\omega\rbrack$, 

\noindent then Res
$\lbrack\omega\rbrack$ is $d$-closed if $\omega$ is $d$-closed.

\vskip 2mm

\noindent {\bf  4. Expression of the residue current of a 
closed meromorphic differential form.}\vskip 2mm

In this section and a part of the following one, we give statements on 
residue currents according to the general hypotheses and proofs of
sections 2 and 3. Proofs in  a particular case where the polar set is
equisingular and the singularity of the polar set is a 2-codimensional
smooth submanifold are given in ($\lbrack$D 57$\rbrack$, IV.D). \vskip 1mm

\noindent {\bf 4.1. Closed meromorphic differential forms. }

\noindent {\bf 4.1.1.} Let $\displaystyle\omega=\frac{\alpha}{f}$ be a
$d$-closed meromorphic differential $p$-form on a small enough open
neighborhood $U$ of  the origin 0 of
$\C^n$. From section 2.1, we get
\ \ \  $\omega= \sum\omega_k$ with \ \ \ $\omega_k=\sum_{\mu=1}^{r_k}{
^j{c_\mu^k}}\frac{\alpha}{{_j\rho_k^\mu}}$ for every $j=1,\ldots,n$. We
have

 $$\displaystyle^j{c_\mu^k}
={\frac{^ja_\mu^k(z_1,\ldots,z_n)}{^jb_\mu^k(z_1,\ldots,\widehat{z_j},
\ldots,z_n)}},$$
where $a$ and $b$ are holomorphic. Then $d\omega=\sum d\omega_k$ and
$d\omega_k$ is the quotient of a  holomorphic form by a product of
$^jb_\mu^k(z_1,\ldots,\widehat{z_j},\ldots,z_n)$ and $_j\rho_k^{r_k+1}$
(see $\lbrack$D 57$\rbrack$, IV,D.1).

As at the end of section 2.2, using the local coordinates
$$(z_1,\ldots,z_{j-1},\rho_k,z_{j+1},\ldots,z_n),$$ we have
$$\omega_k=\sum_{\mu=1}^{r_k}\lbrack _jA_\mu^k\wedge\  _j\rho^{-\mu}_k 
d_j\rho_k +_j\rho^{-\mu}_k B'_k\rbrack, \leqno (4.1)$$
where the coefficients are meromorphic.\vskip 1mm

Let ${\cal R}_j$ be the ring of meromorphic forms on $U$ whose 
coefficients are quotients of holomorphic forms on $U$ by products of 
powers of $\displaystyle\frac{\partial_j\rho_k}{\partial z_j}$ and
$^jb^k_\mu$.

\noindent {\bf Lemma 4.1} ($\lbrack$ D 57$\rbrack$, Lemme 4.10). {\it
Assume that $d\omega_k\in {\cal R}_j$. Then

$$\omega_k= _j\rho_k^{-1}d _j\rho_k\wedge a^k_j+\beta^k_j+dR^k_j$$
 with $$R^k_j=\sum_{\nu=1}^{{r_k}-1} {_je}_\nu^k
{_j\rho}_k^{-\nu}\hskip 2mm and \hskip 2mm da^k_j=d_j{\rho}_k\wedge
\hskip 1mm ^ka_j' +C{_j^k}_j\rho_k,$$ 
where $a^k_j, \beta^k_j,\hskip 1mm _je_\nu^k,\hskip 1mm ^ka_j', C_j \in
{\cal R}_j$ and are independent of $dz_j$.}

\vskip 1mm

\noindent {\bf 4.1.2.} Let $\varphi$ be of type $(n-p,n-1)$. Then 
$$\varphi=\sum\varphi_j, {\rm with}\hskip 2mm \varphi_j=\sum\psi_{l_1,
\ldots,l_{n-p}}dz_{l_1}\wedge\ldots\wedge
dz_{l_{n-p}}\wedge\ldots\wedge\widehat {d\overline
{z}_j}\wedge\ldots$$\vskip 1mm

\noindent {\bf Proposition 4.2.} {\it Let $\displaystyle 
\omega=\frac{\alpha}{f}$ be a $d$-closed meromorphic $p$-form on 
$U$. Given a coordinate system on
$U$, and with notations of section {\rm 2.1}, there exists a current
$S_j^{p-1,1}$ such that $d\sec S_j\vert_{U\setminus Z}=0, \ \ {\rm supp}
S_j=Y$ and, for every $k,j$, a $d$-closed meromorphic $(p-1)$-form
$A^k_j$ on $Y_k$ with polar set $Z$ such that
$$Res\lbrack\omega\rbrack(\varphi)=\sum_{j=1}^n\Big( 2\pi i\sum_k 
Vp_{Y_k,B_j} A^k_j+d'S_j\Big)(\frac{\partial}{\partial z_j} {\rm L}
\varphi_j).$$
When the coordinate system is changed, the first term of the parenthesis 
is modified by addition of 

\noindent $2\pi i \sum_k
d'Vp_{Y_k,B_j}\lbrack F^k_j\rbrack$ where $F^k_j$ is a meromorphic
$(p-2)$-form on $Y_k$ with polar set $Z$}. \vskip 1mm

Here $2\pi i\sum_{j=1}^n\sum_k Vp_{Y_k,B_j} A^k_j(._j)$ will be called
the  {\it reduced residue} of $\omega$.\vskip 1mm

\vskip 1mm

\noindent {\it Proof}. Apply the proof of (*) (section 2) to the 
meromorphic form of Lemma 4.1.

We shall use the expression of ${\rm Res}\lbrack\omega\rbrack(\varphi)$ 
of section 2.2, for $\omega$ closed.

For $k$ and $j$ fixed, we consider 

$$J_{kj}=\lim_{\delta\rightarrow
0}\lim_{\epsilon\rightarrow 0}\int_{\mid B_j\mid\geq\delta, \mid{_j
\rho_k}\mid=\epsilon}\omega_k(\varphi_j).$$
Then $\displaystyle {\rm Res}\lbrack\omega\rbrack
(\varphi)=\sum_{k,j}J_{kj}$.

$\displaystyle\lim_{\delta\rightarrow
0}\lim_{\epsilon\rightarrow 0}\int_{\mid B_j\mid\geq\delta,
\mid{_j\rho_k}\mid=\epsilon}dR_j^k\wedge\varphi_j=(-1)^p\lim_
{\delta\rightarrow
0}\lim_{\epsilon\rightarrow 0}\int_{\mid B_j\mid\geq\delta, \mid{_j
\rho_k}\mid=\epsilon}R_j^k\wedge d\varphi_j.$

Let $S_j^k$ be the current defined by

$$S_j^k(\psi_j)=-\lim_{\delta\rightarrow
0}\lim_{\epsilon\rightarrow 0}\int_{\mid B_j\mid\geq\delta, \mid{_j
\rho_k}\mid=\epsilon}R_j^k\wedge \psi_j.$$

By Lemme 4.1. $R_j^k$ is independent of $dz_j$.

Let $\psi_j= dz_j\wedge\eta^j+\xi^j$, where $\xi^j$ is independent of $dz_j$, then $\eta^j=\frac{\partial}{\partial z_j}{\rm L}\psi_j$.

After change of coordinates:  

$$S_j^k(\psi_j)=-\lim_{\delta\rightarrow
0}\lim_{\epsilon\rightarrow 0}\int_{\mid B_j\mid\geq\delta, \mid{_j\rho_k}
\mid=\epsilon}{\Big(\frac{\partial_j\rho_k}{\partial z_j}\Big
)}^{-1}R_j^k\wedge d_j\rho_k\wedge \eta^j\leqno (4.2)$$

$$=(-1)^p2\pi i\lim_{\delta\rightarrow 0}\sum_\nu\int_{Y_k \vert B_j\vert\geq\delta} (\nu-1)!^{-1}\Big(\frac{\partial^{\nu-1}\big( {_je}_\nu^k\wedge\eta^j\big(\frac{\partial {_j\rho_k}}{\partial z_j}\big)^{-1}\big)}{{\partial {_j\rho}_k}^{\nu-1}}\Big)_{_j\rho_k=0}$$

We have $\displaystyle S_j(\psi_j)=\sum_kS_j^k$.

$\displaystyle\lim_{\delta\rightarrow
0}\lim_{\epsilon\rightarrow 0}\int_{\mid B_j\mid\geq\delta, \mid{_j\rho_k}
\mid=\epsilon}{_j\rho_k}^{-1}d _j\rho_k\wedge a^k_j+\beta^k_j=2\pi i 
\lim_{\delta\rightarrow
0}\int_{\mid B_j\mid\geq\delta} {a^k_j}\mid_ {Y_k}=2\pi i Vp_{Y_k,B_j}
A^k_j$, with $A^k_j={a^k_j}\mid_ {Y_k}$\vskip 1mm

The last alinea is proved as in ($\lbrack$D 57$\rbrack$,
IV.D.4).\hskip 6cm $\square$

\vskip 1mm

\noindent {\bf Corollary 4.3.} {\it The current $S_j$ is obtained by 
application of holomorphic differential operators to currents principal
values of meromorphic forms supported by the irreducible 
components of
$Y$.}\vskip 1mm

\noindent {\it Proof}. The corollary follows from the above expression for
$S_j$ and the computations in section 2.\hskip 2cm $\square$

We remark that $d'$ itself is a holomorphic differential operator.

\vskip 1mm

\noindent {\bf 4.2. Particular cases.}\vskip 1mm

\noindent {\bf 4.2.1. The case} $p=1$. With the notations of Proposition
4.2, the forms
$A^k$ are of degree 0 and are $d$-closed, hence constant and unique: the
reduced residue is  a {\it divisor with complex coefficients}.\vskip 1mm

\noindent {\bf 4.2.2.} With the hypotheses and the notations of
section 2.1, if all the multiplicities $r_k$ are equal to 1, the
reduced residue is uniquely determined and the current $S=0$.\vskip 1mm

\noindent {\bf 4.3. Comparison with the expression of $Res\lbrack\omega\rbrack$ in section 2, when $\omega$ is $d$-closed.}

The reduced residue is equal to
$$2\pi
i\sum_{j=1}^n\Bigl[\sum_k
 Vp\ ^j_{Y_k,B_j}\lbrack
\Big(\frac{\partial\rho_k}{\partial z_j}\Big)^{-1} \ {
^j{c_1^k}}\rbrack\big(\frac{\partial}{\partial z_j}{\rm L}
(\alpha\wedge.)_j\big)\Bigr].$$

It is well defined if all the poles of $\omega$ are simple.

\vskip 2mm     

\noindent {\bf 5. Generalization of a theorem of Picard. Structure of 
residue currents of closed meromorphic forms.}\vskip 2mm  
\noindent {\bf 5.1.}
The theorem of Picard $\lbrack$P 01$\rbrack$  characterizes the divisor 
with complex coefficients associated to a $d$- closed differential form,
of degree 1 of the third kind, on a complex projective algebraic surface;
this result has been generalized by S. Lefschetz (1924): "the divisor has
to be homologous to 0", then by A. Weil (1947). Locally, one of its
assertions is a particular case of the theorem of Dickenstein-Sessa
($\lbrack$DS
 85$\rbrack$, Theorem 7.1): {\it Analytic cycles are locally residual
currents} (see section 5.5), with a variant by D. Boudiaf ($\lbrack$B
92$\rbrack$, Ch.1, sect.3).
\vskip 1mm

\noindent {\bf5.2. Main results}.\vskip 1mm

\noindent {\bf Theorem 5.1.} {\it Let $X$ be a complex manifold which is 
compact K\"ahler or Stein, and $Y$ be a complex hypersurface of $X$, then
$Y=\cup_\nu Y_\nu$ is a locally finite union of irreducible hypersurfaces. Let $Z=$ Sing $Y$, and let $A_\nu$ be a $d$-closed meromorphic $(p-1)$-form on
$Y_\nu$ with polar set $Y_\nu\cap Z$ such that the current $t=2\pi
i\sum_\nu Vp_{Y_\nu} A_\nu$ is $d$-closed.

Then the following two conditions are equivalent:

$(i)$ $t$ is the residue current of a $d$-closed meromorphic $p$-form on $X$ having $Y$ as polar set with multiplicity one.

$(ii)$ $t=dv$ on $X$, where $v$ is a current, i.e., is cohomologous to 0 on
$X$.}\vskip 1mm

\noindent {\it Proof}. From section 4 locally, and a sheaf cohomology 
machinery globally; detailed proof will be given later for the more
general theorem 5.5.\hskip 10cm $\square$\vskip 1mm

For $p=1$, the $A_\nu$ are complex constants, then $t$ is the divisor with complex coefficients $2\pi i\sum_\nu  A_\nu Y_\nu$.\vskip 1mm

\noindent {\bf Corollary 5.1.1.} {\it Under the hypotheses of Theorem
5.1, every residue current of a closed meromorphic p-form appears as a
divisor, homologous to 0, whose coefficients are principal values of
meromorphic
$(p-1)$-forms on the irreducible components of the support of the divisor
and conversely.}\vskip 1mm 

Let ${\cal R}^{loc}_{q,q}(X)$ be the vector space of locally rectifiable 
currents of bidimension $(q,q)$ on the complex manifold $X$ and $${\cal
R}^{loc \C}_{q,q}(X)={\cal R}^{loc}_{q,q}\otimes_{\Z }\C(X)$$.\vskip 1mm

\noindent {\bf Theorem 5.2.} {\it Let $T\in {\cal R}^{loc
\C}_{q,q}(X)$, $dT=0$. Then $T$ is a holomorphic $q$-chain with complex
coefficients.}\vskip 1mm

This is the structure theorem of holomorphic chains of Harvey-Shiffman-Alexander for complex coefficients; thanks to it, divisors will be translated into rectifiable currents.\vskip 1mm

\noindent {\bf Theorem 5.3.} {\it Let $X$ be a Stein manifold or a
compact
 K\"ahler manifold. Then the following conditions are equivalent:  

$(i)$ $T$ is the residue current of a $d$-closed meromorphic $1$-form 
on $X$ 
having {\rm supp} $T$ as polar set with multiplicity $1$;

$(ii)$ $T\in {\cal R}^{loc \C}_{n-1,n-1}(X)$, $T=dV$.}\vskip 1mm

In the same way, we can reformulate the Theorem 5.1 with rectifiable 
currents:
\vskip 1mm

\noindent {\bf Theorem 5.4.} {\it Let $X$ be a Stein manifold or a
compact
 K\"ahler manifold. Then the following conditions are equivalent:  

$(i)$ $T=\sum_\nu a_\nu T_\nu$, with $T_\nu\in{\cal R}^{loc \C}_{n-1,n-1}(X)$, 
$d$-closed, and $a_\nu$ the principal value of a d-closed meromorphic
$(p-1)$-form on {\rm supp} $T_\nu$, such that $T=dV$;

$(ii)$ $T$ is the residue current of a $d$-closed meromorphic $p$-form on $X$ 
having $\cup_l T_l$ as polar set with multiplicity $1$.}

\vskip 1mm
\noindent {\bf 5.3. Remark}. 
 The global Theorem 5.1 gives also local results since any open ball centered at $0$ in  $\C^n$ is  a Stein manifold.\vskip 1mm

\noindent {\bf 5.4. Generalization}.\vskip 1mm

\noindent {\bf 5.4.1.} With the notations of section 4.1, what has been done with the current $2\pi i\sum_\nu Vp_{Y_\nu} A_\nu$ is also possible in the general case. The current $S$ is defined as follows:
 let $\psi=\sum_j \psi_j$, then $S(\psi)=\sum_j\sum_k S_j^k(\psi_j)$.

From (4.2), we have:

$$S_j^k(\psi_j)=2\pi i\sum_{\mu=1}^{r_k}\sum_{l=0}^{\mu-1}
\Delta_{j,k}^{\mu,l}Vp^j_{Y_k,B_j}\lbrack
\gamma^{\mu j}_{k,l}\rbrack\big(\frac{\partial}{\partial z_j}{\rm
L}\psi_j\big)\leqno (5.3)$$ where $\gamma^{\mu j}_{k,l}$ is a meromorphic form on $Y_k$, with polar set contained in $Y_k\cap\{B_j=0\}$,  and where 
$\Delta_{j,k}^{\mu,l}$ is a holomorphic differential operator in the neighborhood of $Y_k$. In the global case, for $Y=\cup_\nu Y_\nu$ locally finite, 
we take $k=\nu$, the sum $\sum_\nu S_j^\nu$ being locally finite.\vskip 1mm

 Then we will get generalizations of the results in sections 5.2 and 5.3 completing the programme of section 1.3.\vskip 1mm
 
 \noindent {\bf Lemma 5.1.} {\it Let $m^p$ be the sheaf of closed
meromorphic differential forms. Let $\overline m^p$ be the image by {\rm
Vp} of $m^p$  in the sheaf of germs of currents on $X$. Then, for 
$X$ Stein  or compact K\"ahler manifold, we have the commutative diagram\vskip 1mm  

$H^0(X,
m^p)\rightarrow H^0(X,\overline
m^p)\rightarrow H^0(X,\overline
m^p/E^p)\rightarrow H^1(X,E^p)$}

$\hskip 4.5cm {\rm Res}\downarrow\hskip 2.5cm  \downarrow$

$\hskip 4.4cmH^0(X,d\sec\overline m^p)\hskip 3mm\rightarrow H^{p+1}(X,\C)$
\vskip 1mm

 (from $\lbrack$D 57$\rbrack$, IV.D.7)\vskip 1mm

 \noindent {\bf 5.4.2.} The residue current of a $d$-closed meromorphic $p$-form is globally written  $t=2\pi i\sum_\nu Vp_{Y_\nu} A_\nu+d'S$,
where $S=\sum_\nu\sum_j S^\nu_j$, with $dt=0$, from the local Proposition  4.2.\vskip 1mm

\noindent {\bf Theorem 5.5.} {\it  If $X$ is a complex manifold which is compact K\"ahler, or  Stein, and $Y$ is a complex hypersurface of $X$, then $Y=\cup_\nu Y_\nu$ is a locally finite union of irreducible hypersurfaces. Let $Z$= {\rm Sing}$Y$; for every $\nu$, let $A_\nu$
be a $d$-closed meromorphic $(p-1)$-form on $Y_\nu$, and, in the notations of $(5.3)$ with $k=\nu$, 
$\gamma^{\mu j}_{\nu,l}$  be meromorphic $(p-2)$-forms on $Y_\nu$, 
with polar set
$Y_\nu\cap Z$ such that the current \ $t=2\pi i\sum_\nu Vp_{Y_\nu} A_\nu +d'S$, \ with $S=\sum_\nu\sum_j S^\nu_j$,  be $d$-closed.

Then the following two conditions are equivalent:

$(i)$ $t$ is the residue current of a $d$-closed meromorphic  $p$-form on $X$ having $Y$ as polar set.

$(ii)$ $t=dv$ on $X$, where $v$ is a current, i.e. $t$ is cohomologous to $0$ on $X$.}\vskip 1mm 

\noindent {\it Proof}. 

$(i)\Rightarrow (ii)$: From Lemma 5.1, the cohomology class of a residue
current is $0$; it is the case of $t$.

\vskip 1mm

$(ii)\Rightarrow (i)$: $t=dv$ on $X$; $t$ of type $(p,1)$ implies:
$t=dv=d\sec v$; $v$ of type $(p,0)$; the current $v$ is closed on 
$X\setminus Y$, therefore it a holomorphic $p$-form on $X\setminus Y$.
Let $m^p_Y$ be  the sheaf of closed meromorphic $p$-forms with polar set
$Y$; the Lemma 5.1 is valid for $m^p_Y$ instead of $m^p$. At a point $O\in
Y$,
$Y$ is defined by $\Pi_k\rho_k=0$ (omitting the index $j$); the $r_k$
being the integers in (5.3), then 
$d(\Pi_k\rho_k^{r_k} v)=\Pi_k\rho_k^{r_k} d\sec v=\Pi_k\rho_k^{r_k} t=0$ from
 Lemma 4.1; therefore $\Pi_k\rho_k^{r_k} v$ is a germ of holomorphic form at $O$ and $v$ extends a closed meromorphic form  $G\in H^0(X,m^p_Y)$ on $X$. 

We will show that $t$ is the residue current of $G$. 

From Proposition 4.2, 
$${\rm Res}\lbrack G\rbrack=d\sec {\rm Vp} \ G=2\pi i\sum_\nu {\rm Vp}_{Y_\nu}B_\nu+d'T$$
where $B_\nu$ and $T$ are of the same nature as $A_\nu$ and $S$.\vskip 1mm

\noindent {\bf Lemma 5.2}. $M=v-$ Vp $G$ \hskip 2mm {\it satisfies \hskip
2mm
$d\sec M=0$}.

\vskip 1mm 
\noindent {\it Proof}. We have:
$$d\sec M= 2\pi i\sum_\nu Vp_{Y_\nu}(A_\nu-B_\nu)+d'(S-T)\leqno (5.4)$$

 Let
$O_1$ be a non singuler  point of $Y$; there exists $k$ such that: $O_1
\in \{{_j\rho}_k=0\}$,
$(j=1,\ldots,n)$; in the neighborhood of  $O_1$, ${_j\rho}_k$ can be used 
as  local coordinate. We have: $M= M_j$  where $M_j$ is
written with the local coordinates
$(\ldots,z_{j-1},{_j\rho}_k,z_{j+1},\ldots)$; $d\sec M=d\sec
M_j$; the support of $d\sec M$ is $Y$, then, in the neighborhood of $O_1$,
$d\sec M_j$ vanishes on the differential forms containing $d{_j\rho}_k$
or $d{_j\overline\rho}_k$. Then
$$d\sec M_j=d{_j\rho}_k\wedge d{_j\overline\rho}_k\wedge N_j \leqno 
(5.5)$$
$M_j$ is of type $(p,0)$, therefore without term in $d{_j\overline\rho}_k$ 
and in $d\overline z_l$, $l\not =j$.

From (5.5), $\displaystyle \frac{\partial M_j}{\partial\overline
z_l}=0$, then $$\displaystyle d\sec
M_j=d{_j\overline\rho}_k\wedge\frac{\partial
M_j}{\partial{_j\overline\rho}_k}\leqno (5.6)$$

$d\sec M_j$ is a differential form with distribution 
coefficients supported by $Y_k$, therefore, outside $Z$,  from the
structure theorem of distributions supported by a submanifold
($\lbrack$Sc 50$\rbrack$, ch. III, th\'eor\`eme XXXVII), and from (5.6), the coefficients of $d\sec M_j$ being those of $\displaystyle \frac{\partial
M_j}{\partial{_j\overline\rho}_k}$, then 
$d\sec M_j$ contains transversal derivatives with respect ${_j\rho}_k$ or 
${_j\overline\rho}_k$ of order at least equal to $r_k+1$, 
what is incompatible with the initial expression (5.4) of $d\sec M_j$,
except if $d\sec M_j=0$ outside $Z$.
From (5.4) the Vp$_{Y_k}(A_\nu-B_\nu)$ and $(S-T)$ being defined as
 limits of integrals of forms vanishing on $Y\setminus Z$, we have: $d\sec
M=0$ on $X$.\hskip 9.5cm$\square$

From Lemma 5.2, Res$\lbrack G\rbrack=d\sec v=t$. \hskip 10cm$\square$
\vskip 1mm
\noindent {\bf Corollary 5.5.1.} {\it Under the hypotheses of Theorem 5.5,
the current $S$ is a sum of currents obtained by application of holomorphic differential operators to principal values of meromorphic forms on the irreducible components $Y_\nu$ of $Y$}.
\vskip 1mm

\noindent {\bf Corollary 5.5.2.} {\it Under the hypotheses of Theorem 5.5,
the residue current of a $d$-closed meromorphic differential $p$-form 
is the sum, cohomologous to 0, of currents obtained by application of
 holomorphic differential operators to currents  
$\sum_\nu a_\nu T_\nu$, with $T_\nu\in{\cal R}^{loc
\C}_{n-1,n-1}(X)$, $d$-closed, and  $a_\nu$ the principal value of a
meromorphic 
 $(p-1)$-form on {\rm supp} $T_\nu$}.
 
 \vskip 1mm
\noindent {\bf 5.5. Remarks}. The Theorems of the sections 5.2 and 5.4 and 
their Corollaries are valid for locally residue currents in the 
terminology of $\lbrack$DS 85$\rbrack$.
Results are also valid for any complex analytic manifold, using less 
natural cohomology (cf $\lbrack$D 57$\rbrack$, IV.D.7).

\vskip 2mm
\noindent {\bf 6. Remarks about residual currents}  
$\lbrack$CH 78$\rbrack$, \noindent $\lbrack$DS 85$\rbrack${\bf .}\vskip 2mm

In the classical definition and notations, we consider residual currents 
$R^p\lbrack\mu\rbrack=R^pP^0\lbrack\mu\rbrack$, where $\mu$ is a
semi-meromorphic  form $\frac{\alpha}{f_1.\ldots.f_p}$, and $\alpha$ a
differential $(p,0)$-form. Then, $R^p\lbrack\mu\rbrack$ satisfies a
formula analogous to (*) of section 2.4. ($\lbrack$D 93$\rbrack$ , section
8).

Locally, one of the assertions of the theorem of Picard is valid for any
$p$, from the result of Dickenstein-Sessa  quoted in section 5.1. So
generalizations of theorems in sections 5.2 to 5.4, for residual
currents, seem valid.\vskip 6mm

\centerline {References}\vskip 2mm
\noindent  [A 97]  H. Alexander, Holomorphic chains and the support hypothesis conjecture, {\it J. of the Amer. Math. Soc.}, {\bf 10} (1997), 123-138.\vskip 1mm
\noindent  $\lbrack$B 92$\rbrack$ D. Boudiaf, {\it Th\`ese de l'Universit\'e Paris VI}, (1992).\vskip 1mm

\noindent  $\lbrack$CH 78$\rbrack$ H. Coleff et M. Herrera, Les courants r\'esiduels associ\'es \`a une forme m\'eromorphe, {\it Springer Lecture Notes in Math.} {\bf 633} (1978).\vskip 1mm

\noindent $\lbrack$DS 85$\rbrack$ A. Dickenstein and C. Sessa, Canonical reprentatives in moderate cohomology, {\it Inv.
Math.} {\bf 80}, 417-434 (1985).\vskip 1mm

\noindent $\lbrack$D 57$\rbrack$ P. Dolbeault, Formes diff\'erentielles et cohomologie sur une vari\'et\'e analytique complexe, II, {\it Ann. of Math.} {\bf 65} (1957), 282-330.\vskip 1mm
\noindent $\lbrack$D 93$\rbrack$ P. Dolbeault, On the structure of residual currents, {\it Several complex
variables, Proceedings of the Mittag-Leffler Institute, 1987-1988}, Princeton Math.
 Notes {\bf 38} (1993), 258-273.\vskip 1mm
 
\noindent  $\lbrack$HS 74$\rbrack$ R. Harvey and B. Shiffman, A characterization of holomorphic chains, {\it Ann. of Math.} {\bf 99} (1974), 553-587.\vskip 1mm

\noindent $\lbrack$HL 71$\rbrack$ M. Herrera and D. Lieberman, Residues and principal values on complex spaces, {\it Math. Ann.} {\bf 194} (1971), 259-294. \vskip 1mm

\noindent  $\lbrack$K 71$\rbrack$ J. King, The currents defined by analytic varieties, {\it Acta Math.} {\bf 127} (1971), 185-220.\vskip 1mm

\noindent $\lbrack$P 01$\rbrack$ E. Picard, Sur les int\'egrales des diff\'erentielles totales de
troisi\`eme esp\`ece dans la th\'eorie des surfaces alg\'ebriques, {\it Ann. Sc. E.N.S.} {\bf 18}  (1901),
397-420.\vskip 1mm

\noindent $\lbrack$Sc 50$\rbrack$, L. Schwartz, Th\'eorie des distributions, new edition, Hermann, Paris 1966.\vskip 1mm

\noindent  $\lbrack$S 83$\rbrack$ B. Shiffman, Complete characterization of holomorphic chains of codimension one, {\it Math. Ann.} {\bf 274} (1986), 233-256.
\vskip 5mm

\noindent Universit\'e Pierre et Marie Curie-Paris 6, I.M.J. (U.M.R. 7586 du C.N.R.S.)\vskip 1mm

\noindent pierre.dolbeault@upmc.fr

\end